	\documentstyle[12pt]{article}
\begin{document}
\title{{ Strong p-completeness of stochastic differential equations
and the existence of smooth flows on noncompact manifolds}}
\author{Xue-Mei Li\thanks{Research supported   by SERC grant GR/H67263.}\\
Mathematics Institute, University of Warwick, Coventry CV4 7AL\\
e-mail: xl@maths.warwick.ac.uk}
\date{}
\maketitle

\newcommand{\A}{{\bf \cal A}}
\newcommand{\B}{{ \bf \cal B }}
\newcommand{\C}{{\cal C}}
\newcommand{\F}{{\cal F}}
\newcommand{\G}{{\cal G}}
\newcommand{\h}{{\cal H}}
\newcommand{\K}{{\cal K}}
\newcommand{\half}{{  {1\over 2}  }}
\newcommand{\heatsemif}{{ {\rm e}^{ \half t\triangle^{h,1}}   }}
\newcommand{\heatsemi}{{ {\rm e}^{\half t \triangle^{h}}   }}

\newtheorem{theorem}{Theorem}[section]
\newtheorem{proposition}[theorem]{Proposition}
\newtheorem{lemma}[theorem]{Lemma}
\newtheorem{corollary}[theorem]{Corollary}
\newtheorem{definition}{Definition}[section]

\def\limsup{\mathop{\overline{\rm lim}}}
\def\liminf{\mathop{\underline{\rm lim}}}
\def\exp{{\rm e}}

 \noindent{\bf   Summary.}
Here we discuss the regularity of solutions of SDE's  and obtain conditions
 under which a SDE on a complete  Riemannian manifold $M$ has a global 
smooth solution flow, in particular improving the usual global Lipschitz
 hypothesis  when $M=R^n$.  There are also results on  non-explosion of
diffusions.


\noindent{\bf Running Title:}  Stochastic Flows.

\section{Introduction}

 Let $M$ be a n-dimensional connected smooth manifold and  $B_t$
 an m-dimensional Brownian motion on a probability space
 $\{\Omega,{\cal F},P\}$  with filtration $\{\F_t\}$.  Consider the 
(Stratonovich) stochastic differential equation(SDE) on $M$: 
	\begin{equation}
dx_t=X(x_t)\circ dB_t+A(x_t)dt.   \label{eq: basic}
\end{equation}

\noindent
 Here $X$ is $C^2$ from $ R^m\times M$ to the tangent bundle $TM$  with $X(x)$:
 $R^m$ $\to T_xM$ a linear map for each $x$ in $M$,  and $A$ is  a $C^2$ vector
 field on $M$. The pair  $(X,A)$ is called a  stochastic dynamical system
(SDS). Let $\{e_1, e_2, \dots, e_m\}$ be an orthonormal basis for $R^m$.
 Set $X^i(x)=X(x)(e^i)$, and write $B_t=(B_t^1,\dots, B_t^m)$. Then 
$(\ref{eq: basic})$ can be written as:
$$dx_t=\sum_{i=1}^m X^i(x_t)\circ dB_t^i +A(x_t)dt.$$
Let  $\{F_t(x)\}$ be the solution to (1) starting from $x$ with explosion time
$\xi(x)$. 

A SDE on a Riemannian manifold  is called a Brownian
 system with drift $Z$ if it  has (i.e. its associated semigroup has) generator
 $\half \triangle +L_Z$. Here $\triangle$ is the Laplacian, $Z$ is a vector 
field  and  $L_Z$ is the Lie derivative in the  direction   $Z$. 
 Its solution is called a Brownian motion with drift $Z$.  Let $h$ be a 
$C^3$ function on $M$.  The Bismut-Witten Laplacian  is
$\triangle^h=:\triangle +2L_{\nabla h}$. 
 A SDE with generator  $\half\triangle^h$ is called a 
h-Brownian system. Its solution is called a h-Brownian motion.

Recall  that a SDE is called {\it  complete} if its
explosion time $\xi(x)=\infty$ for each $x$; it is {\it strongly complete} if
 the solution can be chosen to be  jointly continuous in time and space for
 all  time. Such a  solution is called a {\it continuous flow}.

The known results on the existence of a continuous flow are mostly on $R^n$ 
and  on compact manifolds. On $R^n$ results are given  in terms of
 global  Lipschitz   conditions. See  Blagovescenskii and Friedlin 
\cite{BL-FR}. The  problems concerning the diffeomorphism property of flows 
have been discussed by e.g. Kunita \cite{KUNITA80}, Carverhill and Elworthy 
\cite{CA-EL83}. See  Taniguchi \cite{TANI89} for discussions on the strong
 completeness of a stochastic dynamical system on an open set of $R^n$.  For 
discussions of higher derivatives  of solution flows on $R^n$, see Krylov
 \cite{Krylov} and  Norris \cite{Norr-smc}.

On a compact manifold, a SDE  with $C^2$  coefficients is strongly complete.
  In fact the solution flow is   $C^{r-1}$  if the coefficients are $C^r$.
Moreover the flow consists of diffeomorphisms. See Kunita   \cite{KUNITA80},
Elworthy \cite{EL78}, and Carverhill and Elworthy \cite{CA-EL83}. For 
discussions in the framework of diffeomorphism groups see Baxendale 
\cite{Baxendale80} and Elworthy \cite{ELbook}.

In  the  article, we discuss the regularity of solution flows 
from a new  approach. We introduce the notions of "strong p-completeness". 
Roughly speaking a SDE is strongly p-complete, if the map
$F_\cdot(-)$ is continuous in  time and space for all time  while restricted
 to a smooth $p$-dimensional submanifold of $M$.  This concept   reveals  the
complicated  regularity property of the flow. For example the flow $x+B_t$ on 
$R^n-\{0\}$ is strongly (n-2)-complete but not strongly (n-1)-complete(see  
example 2 in section 2); on $R^n-\ell$, for $\ell$ a smoothly immersed curve
 it is only strongly (n-3)-complete, $n\ge 3$. 

Besides this,  strong 1-completeness turns out to be a powerful tool for
obtaining results on differentiating  semigroups (section 9), for getting
 formulae for the derivatives of the logarithms of the heat kernels
\cite{EL-LI}, or for obtaining related topological and geometrical properties
 of the underlying manifolds\cite{application}\cite{Li.thesis}
 via moment stability. The moment stability part is illustrated in theorem 
\ref{th: p-complete cohomology} below.

\bigskip

\noindent
{\bf Main Results: }

\noindent
{\bf Theorem \ref{pr: n-1 complete}}:
 A stochastic dynamical system on
a smooth manifold  is strongly complete if strongly (n-1)-complete.
\bigskip

Now consider $M$ furnished with a complete Riemannian metric and associated 
Levi-Civita connection.

\noindent{\bf Theorem \ref{th: dim1}/\ref{th: strong completeness} }:
  A SDE on a complete connected Riemannian manifold is strongly p-complete
 if it is complete at one point and its derivative flow $\{T_xF_\cdot\}$
 satisfies: for each compact set $K$ and each $t>0$,
$$\sup_{x\in K}E\left(\sup_{s\le t}|T_xF_s|^{p+\delta}\right)<\infty$$
for some $\delta>0$ ($\delta$ can be taken to be zero for $p=1$).

Note for $p=1$ we only require the first moment of $|T_xF_t|$, so do better
 than a Sobolev  type theorem.

\bigskip

Following from these, we  obtain theorem  \ref{th: existence}  giving criterion
for the existence of a global smooth flow in terms of the coefficients
of the stochastic differential equations. A straightforward
application of theorem \ref{th: existence} extends the standard
global Lipschitz result on $R^n$ (corollary \ref{extension of global}):
denote by $\A$ the differential generator for (1), which is given by
$$\A f(x)=\half \sum_1^m \nabla^2f\left(X^i(x),X^i(x)\right)+A^X(f)(x).$$
Here $A^X=\half \sum_1^m \nabla X^i(X^i)+A$ is the first order part
of the generator.

\noindent{\bf Theorem \ref{pr: dim1}}/{\bf corollary
\ref{extension of global}}:
A  complete SDE on a complete Riemannian manifold is strongly 1-complete 
if 
\[\begin{array}{lll}
H_1(x)(v,v)&=&2<\nabla A^X(v),v>+\sum_1^m<R(X^i,v)(X^i),v>\\
&+&\sum_1^m|\nabla X^i(v)|^2-\sum_1^m{1\over |v|^2}<\nabla X^i(v),v>^2
\end{array}\]
is bounded above. Here  $R$ is the curvature tensor. It is strongly
 complete if $|\nabla X|$ is bounded and if for some constant $c$
$$2<\nabla A^X(v),v>+\sum_1^m<R(X^i,v)(X^i),v>\le c|v|^2.$$

\noindent    There are also more refined results:

\noindent{\bf Theorem \ref{th: 1 in Euclidean space}}
 Let $M=R^n$ with its flat metric. Suppose
the coefficients of the SDE have linear growth, then its solution
flow consists of diffeomorphisms if the first derivatives of its coefficients 
have sub-logarithmic  growth.

\bigskip

Let  $r(x)$ denote the  distance between $x$ and a fixed point in $M$.

\noindent{\bf Theorem \ref{th: h-Brownian}}:
 A Brownian motion with drift $Z$ is complete if the Ricci curvature is bounded
 from below by $-c(1+r^2(x))$, and $dr(Z)\le c(1+r(x))$. It is strongly 
complete if both $|\nabla X|^2$ and $2<\nabla Z(x)(-),-> -Ric_x(-,-)$ have
sub-logarithmic growth in the  distance function $r$.

\bigskip

\noindent{\bf Acknowledgment}
This article developed from  my Warwick university Ph.D. thesis supervised
by D. Elworthy. I would  like to thank him for  suggesting the idea and
and for many useful discussions in particular for improving an earlier
 version of Theorem 2.3  and prof. L. Arnold for helpful comments.  
I am grateful to Prof. ZhanKan Nie  and to
 the Sino-British Friendship Scholarship Scheme for their  support  during my
 Ph.D. work.

\section{Strong p-completeness: definition}

Let  $S_p$ be the space of the  images of all smooth (smooth in the sense of
extending over an open neighbourhood) singular p-simplices. Recall that a 
singular p-simplex in $M$ is a map from the standard p-simplex to $M$. For
 convenience we also use the term singular p-simplex for the image of a 
singular p-simplex map.

 Before giving the definition, here is an example:

\noindent
{\underline {\bf Example 1}} \cite{EL78}, \cite{ELbook} \label{ex: -0 1}
   Let $X(x)(e)=e$, and $A=0$. Consider the following  stochastic differential 
equation  $dx_t=dB_t$ on $R^n$-$\{0\}$ for $n>1$. The solution is: 
$F_t(x)=x+B_t$, which is complete since for a fixed  starting point $x$, 
$F_t(x)$ almost  surely never hits $0$ . But it is not  strongly complete. 
 However for any $n$-$2$ dimensional hyperplane (or a submanifold) $H$ in the 
manifold, $\inf_{x\in H}\xi(x,\omega)=\infty$ a.s., since a Brownian motion does not charge a set of codimension $2$.

\bigskip
This leads to the following definition suggested by D. Elworthy:

\begin{definition} A  SDE on a manifold is  called strongly  $p$-complete 
if its solution can be chosen to be jointly  continuous in time and
 space a.s. for all time when  restricted to a set  $K\in S_p$. 
$\label{def:complete}$
\end{definition}

\noindent
{\underline {\bf Example 2. }}

The example above on $R^n$-$\{0\}$ (for $n>2$) gives us a SDS which is
strongly (n-2)-complete, but not strongly (n-1)-complete.  It is
not strongly (n-1)-complete from proposition ~\ref{pr: n-1 complete}. We
 shall show it is strongly (n-2)-complete.

First note every singular n-2 simplex has an extension to a bounded
Lipschitz map from the cube $[0,1]^{n-2}$ to $M$. Let $U$ be a subset of
$R^{n-2}$ containing a ball radius $\epsilon>0$. Let $f$ be a bounded
Lipschitz map from  $U$ to $R^n$.  We only need to show  that the capacity
Cap(f) of $f(U)$ is zero. For this, the author is   grateful to  Dr P.
 Kr\"oger for the following  proof.    Let $a=\inf_{x\in U}f(x)$.
Clearly Cap(f(U))=0 is equivalent to Cap(2a+f(U))=0.
Thus we may assume $a>0$.
Define $h:R^n\to R\cup \{\infty\}$ as follows:
$$h(y)=\int_U {dx \over |f(x)-y|^{n-2}}.$$
Clearly $h(y)$ is superharmonic. Thus $h(B_t)$ is a supermartingale. By the 
maximal inequality for positive supermartingales, we have:
$$P\{\sup_{0\le s}h(B_s)\ge n\}\le {1\over n} Eh(0).$$
So $P\{\sup_{0\le s}h(B_s)=\infty\}=0.$
This proves $Cap(h^{-1}(\infty))=0$. Next we show $f(U)\subset h^{-1}(\infty)$.
Let $y=f(z)$ for $z\in U$, then for some constant $c$,
\begin{eqnarray*}
h(y)&=&\int_U {dx \over |f(x)-f(z)|^{n-2}}\ge c\int_{U}{dx \over |x-z|^{n-2}}\\
&\ge& c\int_{B_\epsilon} {dx \over |x|^{n-2}}=\infty.
\end{eqnarray*}
Thus Cap(f(U))=0 as wanted. 
\hfill \rule{3mm}{3mm}

\bigskip

For further discussions, we need  the following theorem on the  existence of a
 partial flow, taken from \cite{ELbook} based on  \cite{KUNITA80}. See
also \cite{CA-EL83}.

\begin{theorem} $\label{th: partial flow}$
Suppose $X$ and $A$ are  $C^{r}$, for $r\ge 2$. Then there is a partially
 defined flow $(F_t(\cdot),\xi(\cdot))$ which  is a maximal solution  to 
$(\ref{eq: basic})$ such that if 
$M_t(\omega)=\{x\in M, t<\xi(x,\omega)\}$,
then there is a set $\Omega_0$ of full measure such that for all 
$\omega\in \Omega_0$:
\begin{enumerate}
\item
$M_t(\omega)$ is open in $M$ for each $t>0$, i.e. $\xi(\cdot,\omega)$ is lower
 semicontinuous.
\item
$F_t(\cdot,\omega): M_t(\omega)\to M$ is in $C^{r-1}$ and is a diffeomorphism
 onto an open subset of $M$. Moreover the map : $t\mapsto F_t(\cdot,\omega)$ is
 continuous into $C^{r-1}(M_t(\omega))$, with the topology of uniform convergence
 on compacta of  the first r-1 derivatives.
\item
Let $K$ be a compact set and $\xi^K=\inf_{x\in K} \xi(x)$. Then
\begin{equation}
\lim_{t\nearrow \xi^K(\omega)} \sup_{x\in K} d(x_0, F_t(x))=\infty
\end{equation}
\noindent
almost surely on the set $\{\xi^K<\infty\}$. (Here $x_0$ is a fixed point of
 $M$ and $d$ is any complete  metric on $M$.)
\end{enumerate}
\end{theorem}

From now on, we shall use $(F_t,\xi)$ for the partial flow defined in 
theorem $~\ref{th: partial flow}$ unless otherwise stated. Note that if the
solution can be chosen to be continuous in time and space for all time on
a compact set $K$, then the explosion time $\xi^K$ in the above lemma is
infinite (Elworthy \cite{ELbook}). Thus  strong p-completeness of a SDE
is equivalent to  $\xi^K=\infty$ a.s. for all $K\in S_p$.

\begin{proposition}  If the SDE  considered is  strongly p-complete,
 then  $\xi^N=\infty$ a.s.  for any $p$ dimensional 
smooth submanifold   $N$ of $M$. In particular $F$ can be chosen to be 
$C^{r-1}$ on any given smooth p-dimensional submanifold.
\label{pr: on definition of strong}
\end{proposition}

\noindent
{\bf Proof:}   Let $N$ be a $p$ dimensional submanifold.  Since all smooth 
differential manifolds have a  smooth triangulation (Munkres \cite{Munkres}), 
 we can  write: $N=\cup V_i$. Here $V_i$ are smooth singular p-simplices.  
But  $\xi^{V_i}=\infty$ a.s. for each $i$ from the assumption. Thus 
 $F_\cdot(\cdot)|_{V_i}$ is continuous a.s. and  so is $F|_N$ itself. This
gives  $\xi^N=\infty$  almost surely. The existence of a $C^{r-1}$ 
version comes from a uniqueness result from \cite{ELbook}.
\hfill \rule{3mm}{3mm}

Note that if $\sigma\colon \triangle^p\to M$ is a smooth p-simplex, then by 
\cite{ELbook}, strong p-completeness implies that $F_t\circ\sigma$ has a
 $C^{r-1}$ version.

\bigskip
If $p$ equals the dimension of $M$, strong p-completeness gives back the usual
 definition of strong completeness, i.e. the partial flow defined in theorem 
~\ref{th: partial flow} satisfies $\inf_{x\in M}\xi(x)=\infty$ almost surely. 
In this case we will continue to use strong completeness for strong 
n-completeness.

\begin{theorem}\label{pr: n-1 complete}
A stochastic dynamical system on a n-dimensional manifold is strongly 
complete if   strongly (n-1)-complete. 
\end{theorem}

\noindent
{\bf Proof:} Since we have strong completeness for compact manifolds, we
shall assume $M$ is not compact in the following proof.
Let $B$ be a geodesic ball centered at some point $p$ in $M$ with radius
 smaller 
than the injectivity radius at $p$. Since $M$ can be covered by a countable
 number of such balls, we only need to prove $\xi^B=\infty$ almost surely.

Let $B$ be such a ball. Clearly divides $M-\partial B$ consist of two parts,
one $K_0$ say  bounded and the other $N_0$ unbounded. 
Fix $T>0$. By the ambient isotopy theorem there is a diffeomorphism $H$
from $[0,T]\times M$ to $[0,T]\times M$ given by:
$(t,x)\mapsto (t, h_t(x))$ for $h_t$  some diffeomorphism from $M$ to its
 image, and satisfying:
$$h_t|_{\partial B}=F_t|_{\partial B}.$$
Set $K_t=h_t(K_0)$, $N_t=h_t(N_0)$. Then
$$M=K_t\cup F_t(\partial B) \cup N_t,$$
 and
\begin{equation}
F_t(\stackrel {\circ}{B})\subset K_t
\label{eq: n-complete 1}
\end{equation}
 on $\{\omega: t<\xi^{B}(\omega)\}$.
Now
$$\cup_{0\le t\le T}\bar K_t={Proj}^1\left[H(\bar K_0 \times [0,T])\right],$$
here $\hbox{Proj}^1$ denotes the projection to $M$. Thus 
$\cup_{0\le t\le T} \bar K_t$ is compact. By $(\ref{eq: n-complete 1})$, 
$F_t(B)=F_t(K_0)\cup F_t(\partial B)$, for $0\le t\le T\wedge \xi^B$,  
stays in a  compact region. So $\xi^B\ge T$ almost surely from part 3 of
theorem ~\ref{th: partial flow}.  \hfill \rule{3mm}{3mm}

\bigskip

\noindent{\bf Application of  strong p-completeness}

Let $C^\infty(\Omega^p)$ be the space of $C^\infty$ smooth $p$ forms on $M$,
$H^p(M,R)$  the $p^{th}$  {\it de Rham cohomology} group, and
 $H_K^p(M,R)$ the de Rham cohomolog group for compactly supported  p-forms.
 Recall that a SDS is said to be strongly $p^{th}$-moment stable if
for all $K\subset M$ compact,
$$\mu_K(p)=\limsup_{t\to \infty}\sup_{x\in K}{1\over t}\log E|T_xF_t|^p<0.$$

 The following  theorem  follows from an approach of
\cite{EL-survey} for compact manifolds. For a discussion of such topological
consequences of strong moment stability on noncompact manifolds, see
\cite{application}.

\begin{theorem}\label{th: p-complete cohomology}
Let $M$ be a Riemannian manifold and assume there is a strongly p-complete 
SDS  with strong $p^{th}$-moment stability. Then all bounded closed
 p-forms are exact.  In particular the natural map from $H_K^p(M,R)$ to
 $H^p(M,R)$ is trivial.
 \end{theorem}

\noindent{\bf Proof:}
Let $\sigma$ be a singular p-simplex,  and $\phi$ a bounded closed p-form.
 We shall not distinguish  a singular simplex map from its image.
 Denote by $F_t^*\phi$ the pull back of the form  $\phi$ and 
$(F_t)_*\sigma=F_t\circ\sigma$. Then
$$\int_{(F_t)_*\sigma} \phi =\int_{\sigma} (F_t)^*\phi.$$
But $(F_t)_*\sigma$ is homotopic to $\sigma$ by the strong p-completeness. 
Thus:
$$\int_\sigma\phi=\int_{(F_t)_*\sigma} \phi=\int_\sigma (F_t)^*\phi.$$
Using  strong $p^{th}$ moment stability,
\begin{eqnarray*}
E|\int_\sigma\phi| &=& \lim_{t\to \infty} E|\int_\sigma (F_t)^*\phi| 
\le|\phi|_\infty \lim_{t\to\infty} \,\int_\sigma E|TF_t|^p  \\
&\le&|\phi|_\infty\lim_{t\to\infty} \sup_{x\in \sigma} E|T_xF_t|^p
=0.
\end{eqnarray*}
So $\int_\sigma \phi=0$, and  $\phi$ is exact by  de Rham's theorem.
 \hfill\rule{3mm}{3mm}

Theorem \ref{th: strong completeness} below suggests that strong p-completeness
is not a major restriction given strong moment stability.

\section{Strong 1-completeness}

 Take a sequence of nested relatively 
compact open sets  $\{U_i\}$ such that it is a cover for  $M$ and 
$\bar U_i\subset U_{i+1}$.  Let $\lambda^i$ be a standard  smooth cut off 
function such that:
\[\lambda^i =\left\{ \begin{array} {cl}1  &x\in U_{i+1}  \\
0, & x\not \in U_{i+2}. \end{array} \right.\]
Let $X^i=\lambda^i X$, $A^i=\lambda^i A$, and $F^i_\cdot$ the solution flow to the 
SDS $(X^i,A^i)$. Then $F^i$ can be taken smooth since both $X^i$ and $A^i$ have compact support.  Let $S_i(x)$ 
be the first exit time of $F^i_t(x)$ from $\bar U_i$ and $S_i^K=\inf_{x\in K} S_i(x)$ for a compact set $K$.  Thus 
$S^K_i$ is a stopping time.  Furthermore $F_t^i(x)=F_t(x)$ before  $S_i^K$. 
 Clearly  $S_i^K\le \xi^K$,  and in fact $\lim_{i\to\infty}S_i^K=\xi^K$ as proved in \cite{CA-EL83}.

\bigskip

\noindent  Let 
$$K_1^1=\{\hbox{Image}(\sigma)| \sigma: [0,\ell]\to M \hskip 3pt \hbox{is $C^1$, $\ell<\infty$}\}.$$

Suppose M is given a complete Riemannian metric. Denote by $|-|$ the norm
with respect to this metric. Let $T F_t(v)$ be the derivative of $F_t$
in the direction $v$, whenever it exists. Note it always exists in probability
up to explosion time. See \cite{ELbook}. We shall call $\{TF_t(-): t\ge 0\}$ 
the derivative flow.

\begin{theorem} Let $M$ be a complete connected Riemannian manifold.
Assume there is a point $\bar x\in M$ with $\xi(\bar x)=\infty$ almost surely. Then $\xi^H=\infty$ for all $H\in K_1^1$,  if 

\begin{equation}
\liminf_{j\to \infty}\sup_{x\in K}  E\left(|T_xF_{S^K_j}|\chi_{S^K_j< t}\right)
<\infty
\label{eq: dim1 1}
\end{equation}

\noindent
for every compact set  $K\in K_1^1$ and each $t>0$. In particular when 
$(\ref{eq: dim1 1})$ holds  we have strong 1-completeness, and strong
completeness if the  dimension of  $M$ is less or equal to  $2$.
 \label{th: dim1}
\end{theorem}

\bigskip
\noindent {\bf proof:} 
Let $y_0\in M$.  Let $\sigma_0$ be a piecewise $C^1$ curve parametrized by arc length with end points: $\sigma_0(0)=\bar x$, and $\sigma_0(\ell_0)=y_0$.
Denote by  $K_0$  the image set of the curve.
 Let $K_t=\{F_t(x): x\in K_0\}$, and $\sigma_t=F_t\circ \sigma_0$ be the
composed curve with length $\ell(\sigma_t)$. Then $\sigma_t(\omega)$ is a 
 piecewise $C^1$ curve on $\{\omega: t<\xi^{K_0}(\omega)\}$.
Let $T$ be a stopping time such that  $T<\xi^{K_0}$, then for each $t>0$, 
\begin{eqnarray}
&&E\ell(\sigma_T)\chi_{T<t}
\le E \int_0^{\ell_0} |{d\over ds} \left(F_{T(\omega)}\left (\sigma(s),
\omega\right)\right)|\, ds\, \chi_{T<t}  \\
&&\le \int_0^{\ell_0} E(\chi_{T<t}|T_{\sigma(s)}F_T|)\,ds
\le \ell_0 \sup_{x\in K_0} E\left(|T_xF_T|\chi_{T<t}\right). 
\label{eq: dim1 2}
\end{eqnarray}

\noindent
 Assume $P\{\xi^{K_0}<\infty\}>0$.  There is a number $T_0$ with
 $P\{\xi^{K_0}<T_0\}>0$. On the other hand  there is  also  a
 number $R(\omega)$ such that $R(\omega)<\infty$ a.s.  and
\begin{equation}
\sup_{0\le t\le T_0} d\left(F_t(\bar x,\omega), \bar x\right)\le R(\omega)
\label{flow1} \end{equation}
following from $\xi(\bar x)=\infty$ a.s.  But by theorem 
$\ref{th: partial flow}$, 
\begin{equation}
\lim_{t\nearrow \xi^{K_0}} \sup_{x\in K_0}  d\left(\bar x, F_t(x,\omega)\right)
=\infty 
\label{flow2}  \end{equation} 
 almost surely on  $\{\xi^{K_0}<\infty\}$.
So the triangle inequality combined with $(\ref{flow1})$ and  $(\ref{flow2})$
 yield:
\begin{eqnarray*}
&&\liminf_{t\nearrow \xi^{K_0}}\sup_{x\in K_0} d
\left(F_t(x,\omega),F_t(\bar x,\omega)\right)\ge
\liminf_{t\nearrow \xi^{K_0}} \left[\sup_{x\in K_0} d\left(F_t(x, \omega), 
\bar x \right)-d\left(\bar x, F_t(\bar x, \omega)\right)\right]\\
&&\ge \liminf_{t\nearrow \xi^{K_0}} \sup_{x\in K_0} d
\left(F_t(x,\omega),\bar x\right)
-\sup_{0\le t\le T_0} d\left(\bar x,F_t(\bar x,\omega)\right)
=\infty
\end{eqnarray*}
 on $\{\omega: \xi^{K_0}<T_0\}$. Therefore on this set,
\begin{equation}
\liminf_{t\nearrow\xi^{K_0}}\ell\left(\sigma_t(\omega)\right)
\ge \liminf_{t\nearrow \xi^{K_0}} \sup_{x\in K_0} d\left(F_t(x,\omega),F_t(\bar x,\omega)\right) =\infty
\label{eq: dim1 20}
\end{equation}
 almost surely  for $t\le T_0$.
Let $T_j=: S_j^{K_0}$ be as defined in the beginning of the section, which  converge to $\xi^{K_0}$. Then there is a subsequence, still denoted by 
$\{T_j\}$, s.t. on $\{\xi^{K_0}<T_0\}$,
\begin{equation}
\lim_{j\to \infty} \ell(\sigma_{T_j})\chi_{\xi^{K_0}<T_0}
=\infty, \hskip 4pt a.s.
\label{eq: dim1 15}
\end{equation}

\noindent
However by equation $(\ref{eq: dim1 2})$,  hypothesis $(\ref{eq: dim1 1})$
and Fatou's lemma:
\begin{eqnarray*}
E\liminf_{j\to \infty} \ell(\sigma_{T_j})\chi_{\xi^{K_0}<T_0} &\le&
\liminf_{j\to \infty} E\ell(\sigma_{T_j(\omega)}(\omega))\chi_{T^j<T_0}\\
&\le&   \ell_0 \liminf_{j\to\infty} \sup_{x\in K_0} E|T_xF_{T_j}|\chi_{T_j<T_0}
<\infty,\\  \end{eqnarray*}
contradicting $(\ref{eq: dim1 15})$. 
 Thus $\xi^{K_0}=\infty$. In particular $\xi(y)=\infty$ for all $y\in M$.
 
Next take $K\in K_1^1$, and replace $K_0$ by $K$ in the above proof  to get 
$\xi^K=\infty$,  for  we only used the fact that  there is a point   $\bar x$
in $K_0$ with $\xi(\bar x)=\infty$  and $|TF_t|$ 
satisfies $(\ref{eq: dim1 1})$.  

To see strong 1-completeness, just notice the set of smooth singular
1-simplices $S_1$ is contained in $K_1^1$.
 \hfill \rule{3mm}{3mm}

\bigskip

\noindent
{\bf\underline{Example 3.}} \label{ex: -0 2}

{\bf A.}  The requirement for the manifold to   be  complete is necessary.
e.g.  example   1 on $R^2-\{0\}$ in section 2  satisfies
equation (~\ref{eq: dim1 1}) but is  not strongly complete. In fact if we 
apply the inversion map $z\mapsto {1\over z}$ in complex form  as in \cite{CA-EL83}. The resulting system on $R^2$ is $(\hat X, B)$ where  
\[      \hat X(x,y)=\left[ \begin{array}{cc} y^2-x^2 &2xy\\
-2xy &y^2-x^2\end{array} \right].        \]

The transformed flow  $F_t(z)={z\over 1+zB_t}$  on $R^2$ by inverting does
not satisfy the condition of the theorem on its derivative and it is not
strongly 1-complete.

  \noindent{\bf B.}    Theorem ~\ref{th: dim1} is sharp in the sense it does
 not work if equation (\ref{eq: dim1 1}) is replaced  by 
$\sup_xE|T_xF_t|<\infty$. This can be seen by using the above example
on $M=R^2-\{0\}$ but with the following  Riemannian metric: 
$$|v|^{\#}={|v|\over |x|}, \hskip 6pt v\in T_xM.$$
  This  is  a complete  metric since $\int_0^1 {ds\over s} =\infty$ so 
the point $\{0\}$ is 'at infinity'. But  for each compact set $K$ and  $t>0$
$$\sup_{x\in K}E|T_xF_t|^{\#} =\sup_{x\in K} E{1\over |x+B_t|}
<\infty.$$

\bigskip

We say a SDE is {\it complete at one point}  if there is a point $x_0$ in
$M$ with $\xi(x_0)=\infty$. From the theorem we have the following corollary,
which is known for elliptic diffusions without condition (\ref{eq: dim1 1}).

\begin{corollary}
The SDE (1) is complete if it is complete at one
point and  satisfies condition $(\ref{eq: dim1 1})$ of theorem 
$\ref{th: dim1}$.
\end{corollary}

\section{Strong p-completeness, flows of diffeomorphisms}

\noindent
Denote by $L_p$ the space of all the image sets of  Lipschitz maps
from $[0,1]^p$ to $M$. As in  the last section, we assume that $M$
 is connected and is given a  complete Riemannian metric.

\begin{theorem}\label{th: strong completeness}
Assume  that the SDE (1) is complete at one point. Let $1\le p \le n$.
Then  $\xi^K=\infty$  for each $K\in L_p$, if  for each positive number $t$
and compact set $K$ there is a number  $\delta>0$ such that:
\begin{equation}
\sup_{x\in K} E\left(\sup_{s\le t}|T_xF_s|^{p+\delta}\chi_{s<\xi}\right)
<\infty. 
\label{eq: for strong p-completeness}\end{equation}

\noindent
In particular this implies strong p-completeness.
\end{theorem}

\noindent {\bf Proof:}  Let $\sigma$ be a Lipschitz map from $[0,1]^p$ to $M$
 with image set $K$. Take a compact set  $\hat K $ with the following property:
 for  any two points of $K$, there is a piecewise $C^1$ curve lying in 
$\hat K $ connecting them.

Let $x=\sigma(\underline s)$ and  $y=\sigma(\underline t)$  
and  $\alpha$ be a piecewise $C^1$ curve in $\hat K $ connecting them. 
Denote by $H_\alpha$ the image set of $\alpha$ and $\ell$ its length. 
By proposition $~\ref{th: dim1}$, $\xi^{H_\alpha}=\infty$. Thus for any $T_0>0$
 we have:
\begin{eqnarray*}
&& E\sup_{t\le T_0}\left[d\left(F_{t}(x), F_{t}(y)\right)\right]^{p+\delta}
\le  E\left(\int_0^\ell \sup_{t\le T_0}|T_{\alpha(s)}F_{t}|\, ds\right)^{p+\delta}\\
&\le& \ell^{p+\delta-1}
 E\int_0^\ell \left(\sup_{t\le T_0}|T_{\alpha(s)}F_{t}|^{p+\delta}\right)\, ds
\le  \ell^{p+\delta} 
\sup_{x\in \hat K }\,\left(E\sup_{t \le T_0} |T_xF_{t}|^{p+\delta}\right).
 \end{eqnarray*}

\noindent
Taking infimum over a sequence of such curves which minimizing the distance
between $x$ and $y$,  we get:
$$ E\left(\sup_{t\le T_0}d(F_{t}(x), F_{t}(y))^{p+\delta}\right)
\le d(x,y)^{p+\delta}\sup_{x\in \hat K }\, E\left(\sup_{t\le T_0}
 |T_xF_t|^{p+\delta}\right).$$

The Lipschitz property of the map $\sigma$ gives
$$ E\left(\sup_{t\le T_0}d(F_{t}(\sigma(\underline s),
 F_{t}(\sigma(\underline t)))^{p+\delta}\right)   \le
 c|\underline s-\underline t|^{p+\delta}
\sup_{x\in \hat K }\, E\left(\sup_{t\le T_0} |T_xF_t|^{p+\delta}\right).$$

Thus we have a modification   $\tilde F_\cdot(\sigma(-))$ of
 $F_\cdot(\sigma(-))$ which is jointly  continuous from
 $[0, T_0]\times [0,1]^p\to M$, according to a generalized  Kolmogorov's
 criterion(see e.g. \cite{ELbook}). So for a fixed point $x_0$ in $M$:
$$\sup_{t\in [0, T_0]}\sup_{\underline s\in [0,1]^p} 
d(F_t(\sigma(\underline s), \omega), x_0)<\infty.$$

On the other hand on $\{\xi^K<\infty\}$, 
$\lim_{t\nearrow \xi^K} \sup_{x\in K} d(F_t(x,\omega), x_0)=\infty$
almost surely.  So $\xi^K$ has to be infinity. 

Finally strong p-completeness follows from the fact that every singular
d-simplex has an extension to a Lipschitz map from the cube $[0,1]^p$ to $M$
(by squashing  one half of the cube to the diagonal).
 \hfill \rule{3mm}{3mm}
\bigskip

\noindent{\bf Remarks:}

\noindent
1)  As a consequence, we get that a SDS is strongly complete if it is
 complete at one point and satisfies:
$$\sup_{x\in K}E\sup_{s\le t}|T_xF_s|^{n-1+\delta}<\infty$$
for some $\delta>0$ and for each compact subset $K$ of M. On the other hand,
any direct application of a Sobolev type inequality would require that the
above integrability condition  holds  for a  $p^{th}$  power ($p>n$) of
$|T_xF_t|$.

\noindent
2) Note condition (\ref{eq: for strong p-completeness}) in the theorem cannot 
be replaced by $\sup_{x\in K}E|T_xF_t|^{p+\delta}$  is finite, since the flow 
$x+B_t$ with the  complete Riemannian metric  $<,>^\#$
 in example 3 (section 3) satisfies: for $p<n$,
$\sup_{x\in K}E\left(|T_xF_t|^\#\right)^p<\infty$.

\bigskip

\noindent{\bf Flows of diffeomorphisms}

For the diffeomorphism  property, we only need to look at the ``adjoint''
 system of (\ref{eq: basic}):
\begin{equation}  dy_t=X(y_t)\circ dB_t -A(y_t)dt.
\label{eq: adjoint}
 \end{equation}
A  strongly complete SDE has a flow of diffeomorphisms if and only its
 adjoint equation is also strongly complete. See Kunita \cite{KUNITA80}. See
 also  Carverhill and Elworthy \cite{CA-EL83}.
 Suppose there is a uniform cover for $(X,A)$. Then  its flow consists of 
diffeomorphisms  if for each compact set $K$,
$$\sup_{x\in K}E\sup_{s\le t} \left(|T_xF_s|^{n-1+\delta} 
+\left(|T_{F_s^{-1}(x)}F_s|^{-1}\right)^{n-1+\delta} \right)<\infty,$$
since  in this case both  equation (\ref{eq: basic}) and (\ref{eq: adjoint}) 
are strongly complete by the previous theorem. In this case we also have 
the $C_0$-property, i.e. the associated semigroup preseves $C_0(M)$, the space
 of continuous functions vanishing at infinity. See \cite{EL82}.

\section{Existence of smooth flows}

 Let $M$ be a Riemannian manifold with Levi-Civita connection $\nabla$. There 
is the   stochastic covariant differential equation 
\begin{equation}
dv_t=\nabla X(v_t)\circ dB_t +\nabla A(v_t)dt.
\label{eq: covariant}
\end{equation}

Denote by  $T_{x}F_t(v)$ its solution starting from $v$.  It is in fact the
derivative of $F_t(x)$ in measure.  See \cite{ELbook}.
Let $x_0\in M$,   $v_0\in T_{x_0}M$.  We shall write  $x_t=F_t(x_0)$, and
$v_t=T_{x_0}F_t(v_0)$.

The expectations of the norms of $|v_t|$ can be estimated through the following
equation (see e.g.  Elworthy\cite{ELflow}, or\cite{ELflour}):
\begin{equation}
\begin{array}{ll}
|v_t|^p=|v_0|^p &+p\sum_{i=1}^m\int_0^t|v_s|^{p-2} <\nabla X^i(v_s), v_s> dB_s^i\\
&+{p\over 2}\int_0^{t}|v_s|^{p-2}H_p(x_s)(v_s, v_s)ds.
\end{array}
\label{eq: vt}
\end{equation}
on $\{t<\xi\}$.  Here 
\begin{equation}
 \begin{array}{ll}
H_p(x)(v,v)&= 2<\nabla A(x)(v), v> +\sum_{i=1}^m <\nabla^2 X^i(X^i,v), v>\\
&+\sum_1^m <\nabla X^i(\nabla X^i(v)),v> +\sum_1^m |\nabla X^i(v)|^2\\
&+(p-2)\sum_1^m {1\over |v|^2} <\nabla X^i(v), v>^2,
\end{array}
\end{equation}
for all $x\in M$ and  $v\in T_xM$.  To simplify notation, let
\begin{eqnarray}
M_t^p&=& \sum_1^m p\int_0^t  {<\nabla X^i(v_s), v_s> \over |v_s|^2} dB_s^i, \\
a_t^p&=& {p\over 2} \int_0^t {H_p(x_s)(v_s,v_s) \over |v_s|^2} ds. 
\end{eqnarray}
Here $M_t^p$ and $a_t^p$ depends on the point $(x_0,v_0)\in TM$. We shall 
omit the superscript $p$ if  no confusion is  caused.
Then equation (\ref{eq: vt}) gives:
\begin{equation}\label{eq: vt as exponential}
|v_t|^p=|v_0|^p\exp^{M_t^p-{<M^p,M^p>_t \over 2}+a_t^p}
\end{equation}
 as used by Taniguchi\cite{TANI89}.  Let $|X(x)|^2=\sum_1^m|X^i(x)|^2$, and let
$|\nabla X(x)|^2=\sum_1^m |\nabla X^i(x)|^2$.   We have:

\begin{theorem}\label{th: existence}
Let $M$ be a complete connected Riemannian manifold. Suppose the SDE (1) 
is complete at one point. Let $p>0$. Assume there is a function 
$f: M\to [0,\infty)$ such that:
\begin{enumerate}
\item
$\sup_{x\in K}E\left(\exp^{6p^2 \int_0^t f(F_s(x))\chi_{s<\xi(x)}ds}\right)<\infty$, 
for all $t>0$, $K$ compact.
\item
$|\nabla X(x)|^2\le f(x)$.
\item
 $H_p(x)(v,v)\le 6pf(x)|v|^2$ for all $x\in M$ and $v\in T_xM$. 
\end{enumerate}
Then the system is complete and
 $$E\left(\sup_{s\le t}|T_xF_s|^p\right)<c
E\left(\exp^{6p^2 \int_0^{t} f(F_s(x)) ds}\right).$$
In particular  the system is strongly  d-complete for $d<p$.
\end{theorem}

\noindent{\bf Proof:} First we assume  that the SDE is complete.
Applying Schwartz's inequality to
 equation (\ref{eq: vt as exponential}), we get for each $p>0$:
\begin{eqnarray*}
E\left(\sup_{s\le t} |v_s|^p\right)
\le |v_0|^p\left(E\sup_{s\le t}\exp^{2M_s-<M,M>_s}\right)^\half
\left( E\sup_{s\le t}\exp^{2a_s}\right)^\half.
\end{eqnarray*}

\noindent
 Since $$E\left(\exp^{6 <M,M>_s}\right)
\le  E\left(\exp^{6 p^2 \int_0^t f(x_s)ds}\right)<\infty,$$
 $\exp^{2M_s -{<M,M>_s\over 2}}$ is a martingale by Novikov's criterion 
\cite{RE-Yor91}. Consequently
\begin{eqnarray*}
&&E\left(\sup_{s\le t} \exp^{2M_s -<M,M>_s}\right)\le
4\sup_{s\le t} E\left(\exp^{2M_s-<M,M>_s}\right)\\
&=& 4\sup_{s\le t} E\left(\exp^{2M_s-4<M,M>_s} \exp^{3<M,M>_s}\right)
\le 4\left[E\left(\exp^{6<M,M>_t}\right)\right]^\half,
\end{eqnarray*}
by Cauchy Schwartz and using the fact that  $\exp^{4M_s-8<M,M>_s}$ is a
 supermartingale.  Also
\begin{eqnarray*}
E\left(\sup_{s\le t}\exp^{2a_s}\right)=
 E\left(\sup_{ \alpha\le t}\exp^{p\int_0^\alpha
 {H_p(v_s,v_s)\over |v_s|^2} ds}\right)
\le E\left(\exp^{6p^2\int_0^t f(x_s)ds}\right)
\end{eqnarray*}
giving
$$E\left(\sup_{s\le t}|v_s|^p\right) 
\le 2|v_0|^p\left[E\left(\exp^{6p^2\int_0^tf(x_s)ds}\right)\right]^{3\over 4}
<\infty.$$

Thus for some constant $c_2$(depending only on p and n),
\begin{eqnarray*}
E\left(\sup_{s\le t} |T_xF_s|^p\right)
\le c_2 E\left(\exp^{\left(6p^2 \int_0^t f(F_s(x)) ds \right)}\right).
\end{eqnarray*}

  Thus for each compact subset $K$ of the manifold,
$$\sup_{x\in K} E\left(\sup_{s\le t} |T_xF_s|^p\right)
\le c_2\sup_{x\in K} E\left(\exp^{\left(6p^2 \int_0^t f(F_s(x)) ds
 \right)}\right) <\infty. $$

\bigskip

Next assume (1) is complete at one point, we shall show that it is complete
 everywhere. 
Let $K\in K_1^1$ be a compact subset of $M$,  and let $S_j^K$ be stopping times
 as in  theorem ~\ref{th: dim1}. Then
$$|T_xF_{t\wedge S_j^K}(v_0)|=|v_0|\exp^{\left(M_{t\wedge S_j^K}
-{<M,M>_{t\wedge S_j^K}\over 2}
+a_{t\wedge S_j^K}\right)}.$$

\noindent
Similar calculations as above yield:
\begin{eqnarray*}
\sup_{x\in K}E\left(|T_xF_{S_j^K}|\chi_{S_j^K<t}\right)
&\le& c\sup_{x\in K} E\left(\exp^{6p^2 \int_0^{S_j^K} 
f(F_s(x)) ds}\chi_{S_j^K<t}\right) \\
&\le& c \sup_{x\in K} E\left(\exp^{6p^2 \int_0^{t} 
f(F_s(x)) ds}\right)<\infty.
\end{eqnarray*}
Here $c$ is a constant. The completeness follows from theorem ~\ref{th: dim1}. The strong completeness 
follows from theorem ~\ref {th: strong completeness}.
\hfill \rule{3mm}{3mm}

\bigskip
 
Note that the first condition in theorem ~\ref{th: existence} is a workable
 condition, since Jensen's inequality gives:
\begin{equation}
\label{eq: Jensen}
 E\exp^{\left(6p^2\int_0^t f^2(x_s)\chi_{s<\xi}ds\right)}
\le {1\over t}\int_0^t E\left[\exp^{6p^2tf(x_s)\chi_{s<\xi}}\right]ds. 
\end{equation}
  
For example, take $f\equiv 1$ in  theorem  ~\ref{th: existence}. Let
 $A^X=\half \sum_1^m \nabla X^i(X^i)+A$, and let $R$ be the curvature tensor on
$M$. Recall that the differential generator $\A$ is given by
$$\A f(x)=\half \sum_1^m \nabla^2f\left(X^i(x), X^i(x)\right)
+A^X(f)(x).$$  We next see that  the theorem is a direct extension of the
 global Lipschitz results for $R^n$:

\begin{corollary}\label{extension of global}
The  SDE (1) is strongly complete if it is
complete at one point and  satisfies:  $|\nabla X|$ is bounded and 

 $$2<\nabla A^X(v),v>+\sum_1^m<R(X^i,v)(X^i),v>\le c|v|^2$$
for some  constant $c$.  In fact under these conditions,
$$\sup_{x\in M}E\left(\sup_{s\le t}|T_xF_s|^p\right)<\infty, \hskip 4pt
\hbox{ for all} \hskip3pt p.$$
 The solution to (\ref{eq: basic}) consists of 
diffeomorphisms if also  the ``adjoint'' equation (\ref{eq: adjoint})
 is complete at one point  and 
$$2<\nabla (-A)^X(v),v>+\sum_1^m<R(X^i,v)(X^i),v>\le c|v|^2.$$
\end{corollary}

\noindent{\bf Proof}: \hskip 4pt  First  
$$\nabla A^X(v)=\half \sum_1^m\nabla^2 X^i(v,X^i)
+\half \sum_1^m \nabla X^i(\nabla X^i(v))+\nabla A(v).$$
But by definition of the curvature tensor,
$$<\nabla ^2 X^i(X^i,v),v>-<\nabla^2X^i(v,X^i),v>
=<R(X^i,v)(X^i),v>.$$

\noindent
So
\[\begin{array}{lll}
<\nabla A^X(v),v>&=&<\nabla A(x)(v),v>-\half \sum_1^m<R(X^i,v)(X^i),v>\\
&+&\half \sum_1^m <\nabla^2 X^i(X^i,v),v>+
\half \sum_1^m <\nabla X^i(\nabla X^i(v)),v>.
\end{array}\]

\noindent
Thus
\[\begin{array}{lll}
H_p(x)(v,v)&=&2<\nabla A^X(v),v>+\sum_1^m<R(X^i,v)(X^i),v>\\
&+& \sum_1^m |\nabla X^i(v)|^2+(p-2)\sum_1^m {1\over |v|^2}<\nabla X^i(v),v>^2.
\end{array}\]

Note the last two terms of $H_p$ are bounded. Thus the conditions of theorem 
\ref{th: existence} are satisfied, and the SDE is strongly complete.
 For the diffeomorphism property, note that  the 'adjoint' equation  has
\[\begin{array}{lll}
H_p(x)(v,v)&=&2<\nabla (-A)^X(v),v>+\sum_1^m<R(X^i,v)(X^i),v>\\
&+& \sum_1^m |\nabla X^i(v)|^2+(p-2)\sum_1^m {1\over |v|^2}<\nabla X^i(v),v>^2,
\end{array}\]
and is thus also strongly complete. 
\hfill\rule{3mm}{3mm}
\bigskip

\noindent
However for strong 1-completeness, we can do better:

\begin{theorem}\label{pr: dim1}
Assume (1) is complete at one point, and $H_1(x)(v,v)\le c|v|^2$ for some
 constant $c$. Then we have strong 1-completeness for (1). Furthermore if
 the dimension of $M=2$, then it is strongly complete.
\end{theorem}

\noindent
{\bf Proof:} Let $K\in K_1^1$, and $S_j^K$ be the  corresponding stopping times
as in theorem  ~\ref{th: dim1}. Then 
\begin{equation}
\begin{array}{ll}
|v_{t\wedge S_j^K}|=|v_0| &+\sum_{i=1}^m\int_0^{t\wedge S_j^K}|v_s|^{-1} <\nabla X^i(v_s), v_s> dB_s^i\\
&+{1\over 2}\int_0^{t\wedge S_j^K}|v_s|^{-1}H_1(v_s, v_s)ds.
\end{array}
\end{equation}
from equation (\ref{eq: vt}) with $t$ replaced by $t\wedge S_j^K$, and letting $p=1$. 
On the other hand, 
$$|T_xF_{t\wedge S_j^K}(v_0)|=|v_0|\exp^{\left(M_{t\wedge S_j^K}-
{<M,M>_{t\wedge S_j^K}\over 2} +a_{t\wedge S_j^K}\right)},$$
by (\ref{eq: vt as exponential}).
But $<M,M>_{t\wedge S_j^K}$ and $ a_{t\wedge S_j^K}$ are both bounded, since both
$|\nabla X^i(x)|$ and $H_1(x)$  are bounded on compact sets.
So $|T_xF_{t\wedge S_j^K}(v)|$ is bounded for each $j$ and $v\in T_xM$.
Thus
$$E\int_0^{t\wedge S_j^K} |v_s|^{-1}<\nabla X^i(v_s), v_s> dB_s^i=0.$$
Therefore,
\begin{eqnarray*}
E|T_xF_{t\wedge S_j^K}(v_0)|&=&|v_0| 
+{1\over 2} E\int_0^{t\wedge S_j^K}|v_s|^{-1}H_1(v_s, v_s)ds \\
&\le& |v_0|+  \half c\int_0^t  E|T_xF_{s\wedge S_j^K}(v_0)| ds
\end{eqnarray*}

\noindent  Gronwall's inequality gives:
$ E|T_xF_{S_j^K\wedge t}(v_0)| \le |v_0|\exp^{ct/2}.$
So
\begin{equation}
E\left(|T_xF_{S_j^K}|\chi_{S_j^K<t}\right)\le E|T_xF_{S_j^K\wedge t}|
 \le \exp^{ct/2}.
\label{stopping estimate}
\end{equation}
The strong 1-completeness follows from theorem ~\ref{th: dim1},  and the strong 
completeness for 2-dimensional manifolds follows from theorem ~\ref{pr: n-1 complete}
\hfill\rule{3mm}{3mm}

\bigskip

It is possible to get a slightly different result from theorem 
\ref{th: existence} using the fact that
$$|v_t|^p=|v_0|^p\exp^{M_t^p}\exp^{{p\over 2}\int_0^t{\tilde H(x_s)(v_s,v_s)
\over |v_s|^2}},$$
where 
\begin{equation}
 \begin{array}{ll}
\tilde H(x)(v,v)&= 2<\nabla A^X(x)(v), v> +\sum_{i=1}^m <R(X^i,v)(X^i),v>\\
&+\sum_1^m |\nabla X^i(v)|^2
-2\sum_1^m {1\over |v|^2} <\nabla X^i(v), v>^2,
\end{array}
\end{equation}
and the fact \cite{RE-Yor91}
 $$E\sup_{s\le t}\exp^{\alpha M_s^p}\le
E\exp^{\sup_{s\le t}\alpha M_s^p}<\infty,$$
if $E\exp^{2 \alpha^2 <M^p,M^p>_t}<\infty$. So just as before, if
$|\nabla X|$ is bounded, then we have strong completeness if $\tilde H$
is bounded above. This allows consequent variations in the results below.

\section{Existence of flows on $R^n$}

In this section we shall show some direct consequences of theorem 
~\ref{th: existence}. The  usual global Lipschitz condition 
is improved to allow some growth of the derivatives of the coefficients
(see theorem \ref{th: 1 in Euclidean space}). Consider on $R^n$
\begin{equation}\label{eq: basic on E}
\hbox{(It\^o)} \hskip 20pt  dx_t=X(x_t)dB_t + A(x_t)dt.
\end{equation}

\noindent
It can be rewritten in Stratonovich form:
$$dx_t=X(x_t)\circ dB_t +\bar A(x_t)dt,$$
where $\bar A=A-\half \sum_1^m D X^i(X^i)$.
 So
\begin{equation}\label{H for Euclidean}
\begin{array}{ll}
H_p(v,v)=2<D A(v),v> + |D X(v)|^2
 +(p-2) \sum_1^m {1\over |v|^2} <D X^i(v), v>^2.
\end{array}\end{equation}
Thus the second derivative of $X$ is not involved. 
Let $g\colon R^n \to [0,\infty)$ be a $C^2$
function. Then by It\^o's formula, on $\{t<\xi\}$
\begin{equation}\label{250}
\exp^{g(x_t)}= \exp^{g(x_0)+N_t-{<N,N>_t \over 2} +b_t},
\end{equation}
where $N_t=\int_0^t Dg(X(x_s)dB_s)$ and
\[\begin{array}{lll}
b_t&=& \int_0^t  \half\sum_1^m \left(\left[(Dg)(x_s)(X^i(x_s))\right]^2 +
 (D^2g)(x_s)(X^i(x_s),X^i(x_s))\right)ds \\
&+&\int_0^t (Dg)(x_s)\left(A(x_s)\right)ds.
\end{array}\]

\begin{lemma} Let $c$ be a constant. Let $\tau$ be a stopping time
 with $\tau<\xi$ on $\{\xi<\infty\}$. Then for some constant $k$
$$E\exp^{\left(cg(x_{t\wedge\tau})\right)} 
\le \exp^{c\left(g(x_0)+kt\right)},$$
 provided that 
 $$\half \sum_1^m|Dg(X^i)|^2+\half \sum_1^mD^2g(X^i,X^i)+Dg(A)
\hskip 16pt \hbox{is bounded above.}  $$ 
\end{lemma}

\noindent {\bf Proof:}
 Replacing $t$ by $t\wedge \tau$ in (\ref{250}),
and $g$ by $cg$, then taking expectations on both sides of the inequality
 above,  we get the  required inequality.

\begin{theorem} \label{th: 1 in Euclidean space}
 The SDE  (\ref{eq: basic on E}) on $R^n$  with $C^2$ coefficients is strongly 
complete if  its coefficients have linear growth(in an  extended sense), i.e.
\begin{eqnarray*}
|X(x)|&\le& c(1+|x|^2)^{1\over 2}\\
<x, A(x)>&\le& c(1+|x|^2),
\end{eqnarray*}
and the derivatives of the coefficients have sub-logarithmic growth, i.e.
\begin{eqnarray}
|\nabla X(x)|^2 &\le& c[1+\ln(1+|x|^2)]\\
<\nabla A(x)(v),v> &\le& c[1+\ln(1+|x|^2)]|v|^2
\end{eqnarray}
for all $x$ and $v\in R^n$. Here  $c$ is a constant. In fact under
these conditions  we have:
$$E|x_t|^{2p} \le c_{1,p}\left(1+|x_0|^2\right)^p
\exp^{c_{2,p}t}$$
for some constant $c_{1,p}$ and $c_{2,p}$  depending only on $p$ and
$\sup_{x\in K}E\sup_{s\le t}|T_xF_s|^p$ is finite  for all  $p$ and 
 compact sets $K$.
\end{theorem}

\noindent
{\bf Proof:}
Let $f(x)=[1+\ln(1+|x|^2)]$, $g(x)=\ln(1+|x|^2)$. Then
$$Df(x)\left(A(x)\right)=Dg(x)\left(A(x)\right)
={2<x,A(x)>\over 1+|x|^2},$$
and
\[\begin{array}{lll}
&&D^2f(x)\left(X^i(x),X^i(x)\right)=
{2<X^i(x),X^i(x)>\over 1+|x|^2} -{4<x,X^i(x)>^2\over (1+|x|^2)^2}.
\end{array}\]
So  by the previous lemma (applied to  the function $g$),
$$E|x_{t\wedge T}|^2\le (1+|x_0|^2)\exp^{k_1t}-1$$
for some constant $k_1$ and stopping times $T$ with $T<\xi$. Thus the system 
is complete by a standard argument.
 Applying the same lemma  to $cf$, we have:
$$E\exp^{c[1+\ln(1+|x_t|^2)]} \le \exp^c(1+|x_0|^2)^c\exp^{kt}$$
for some constant $k$($k$ may depend on $c$).  So
$$\sup_{x\in K} E\left(\exp^{6p^2 \int_0^t c[1+\ln(1+|x_s|^2)]ds}\right)
=\sup_{x\in K} {1\over t} \int_0^t \exp^{6p^2ct}(1+|x_0|^2)\exp^{ks}ds
<\infty.$$
 The strong completeness  follows  from theorem ~\ref{th: existence}, using
(\ref{H for Euclidean}) and the assumptions on $\nabla X$ and $\nabla A$.
\hfill \rule{3mm}{3mm} 

\bigskip

For related estimates on $E|x_t|^p$, see \cite{Krylov}. Note  that there is a
stochastically complete SDE on $R^2$ with $|\nabla X(x)|\le |x|$ but which
is  not strongly complete: let $A\equiv 0$, and
$X(x,y)=\left(\begin{array}{cc}y&0\\0&{x^2\over 2}\end{array}\right)$.
  See  Kunita\cite{Kunitabook}.

A different choice of the function $f$ in theorem \ref{th: existence} leads 
to an improvement of a theorem of Taniguchi \cite{TANI89}:

\begin{corollary}
The SDE (\ref{eq: basic on E}) on $R^n$ is strongly complete if for some
 $\epsilon\ge 0$:
\[\begin{array}{lll}
|X^i(x)|&\le& c(1+|x|^2)^{\half -\epsilon}\\
<x, A(x)> &\le& c (1+|x|^2)^{1-\epsilon}
\end{array}\]
\[\begin{array}{lll}
|DX^i(x)|^2&\le& c(1+|x|^2)^\epsilon\\
<\nabla A(x)(v),v>&\le&c (1+|x|^2)^\epsilon|v|^2.
\end{array}\]
\end{corollary}

\noindent{\bf Proof:}  \,\, Clearly the stochastic differential equation
 is complete.
Take $g(x)=c(1+|x|^2)^\epsilon$ in the lemma for $\epsilon>0
(\epsilon=0$ gives the usual globally Lipschitz  continuous condition). Then
\[\begin{array}{lll}
Dg(x)\left(A(x)\right)&=&2c\epsilon(1+|x|^2)^{\epsilon-1}<x, A(x)>,\\
D^2g(x)\left(X^i(x),X^i(x)\right)&=&2c\epsilon (1+|x|^2)^{\epsilon-1}<X^i(x),X^i(x)>\\
 &&+4c\epsilon(\epsilon-1)(1+|x|^2)^{\epsilon-2}<x, X^i(x)>^2.\\
\end{array}\]
So  lemma 6.1 applies to get   $E\exp^{cg(x_t)}<\exp^{c(g(x_0))+2kt}$
 for some constant $k$ and the result follows from theorem 
\ref{th: existence}.  
  \hfill\rule{3mm}{3mm}

This theorem improves a theorem of Taniguchi since: 
(a)  We only need growth conditions on the normal parts of $A$ and $\nabla A$, and
(b)  we do not assume $\epsilon>{1\over 3}$ as in \cite{TANI89}.

\section{Existence of  flows on manifolds with a pole}

 A similar argument  on the existence of flow (c.f. theorem ~\ref
{th: 1 in Euclidean space}) to that on $R^n$ works for  general manifolds  to 
allow the coefficients to have  unbounded derivatives. We first assume  that 
$M$ is  equipped with a Riemannian metric such that there is a pole $P$ in $M$,
i.e.  the distance function $r(-)\colon M \to R$ from $P$ is  smooth. Recall
that $A^X=\half \sum_1^m\nabla X^i(X^i)+A$.

\begin{theorem}\label{th: pole}
Let $M$ be a complete Riemannian manifold with a pole. 
Assume the sectional curvature is bounded from below by $-L^2(r(-))$.  Here
$L$ is a nondecreasing function bigger or equal to $1$. Then the SDE (1.1) 
is complete   and
$$E\left[r(x_t)\right]^p\le \left[1+r(x_0)\right]^p \exp^{k_0[1+p^2]t}$$
 for some constant $k_0$, if the following holds for some constant $c$:
\begin{enumerate}
\item
$|X(x)|^2\le {c[1+r(x)]\over L(r(x))\hbox {coth}\left(r(x)L(r(x))\right)}$;
\item
 $dr(A^X(x))  \le { c}[1+r(x)]$.
\end{enumerate}
\noindent
 It is strongly complete  and
 $\sup_{x\in K}E\sup_{s\le t}|TF_s|^p<\infty$
 for all $p$ and compact sets $K$, if  we also have:

3.  $|\nabla X(x)|^2\le c[1+\ln(1+r(x))]$;

4. $2<\nabla A^X(v),v>+\sum_1^m <R(X^i,v)(X^i),v>\le c[1+\ln(1+r(x))]|v|^2$.

\end{theorem}

\noindent
{\bf Proof:}   \hskip 4pt  First we have:
\[\begin{array}{lll}
r(x_t)&=&r(x_0)+\int_0^t dr(X(x_s)dB_s)
+\half \sum_1^m \int_0^t\nabla^2r(X^i(x_s),X^i(x_s))ds\\
&&+\int_0^t dr(A^X(x_s))ds. \end{array}\]

\noindent
But by Hessian comparison theorem in \cite{Greene-Wu}(p.19 and 
example 2.25 on p.34. The results there is for constant $L$, but
the proof depends only on the behaviour of the manifold around
the geodesic from $p$ to $x$),
$$\nabla^2r(x)\le L(r(x)) \hbox{coth}\left(r(x)L(r(x))\right).$$

Let $T_n(x)$ be the first exit time of $F_t(x)$ from the geodesic ball 
$B(p,n)$, centered at $p$ and radius $n$. Then
\begin{eqnarray*}
Er(x_{t\wedge T_n})&=& r(x_0)
+\half \sum_1^m E\int_0^{t\wedge T_n}\nabla^2r(X^i(x_s),X^i(x_s))ds\\
&&+E\int_0^{t\wedge T_n}dr(A^X(x_s))ds\\
&\le& r(x_0)+{k_1\over 2}\int_0^t E\chi_{s<T_n}(1+r(x_s))ds.\\
\end{eqnarray*}
Here $k_1$ is a constant. Thus
$$Er(x_{t\wedge T_n})\le [r(x_0)+{k_1t\over 2}]\exp^{k_1t/2}.$$

So
\begin{eqnarray*}
&&P\{T_n<t\}={1\over n} E\left(r(x_{t\wedge T_n})\chi_{T_n<t}\right)\\
&\le&{1\over n}[r(x_0)+{k_1t\over 2}]\exp^{k_1t/2}\to 0,\\
\end{eqnarray*}
as $n$ goes to infinity. Thus there is no explosion. Now
 \begin{eqnarray*}
[1+r(x_t)]^p&=&[1+r(x_0)]^p+p\int_0^t [1+r(x_s)]^{p-1}dr(X(x_s)dB_s)\\
&+&{p(p-1)\over 2}\sum_1^m \int_0^t [1+r(x_s)]^{p-2} [dr(X^i(x_s)]^2ds \\
&+&{p\over 2} \sum_1^m \int_0^t [1+r(x_s)]^{p-1}
\nabla^2r(X^i(x_s),X^i(x_s))ds\\
 &+& p\int_0^t [1+r(x_s)]^{p-1}dr(A^X(x_s))ds. \end{eqnarray*}

\noindent
Let $$M_t=\int_0^t p{dr(X(x_s)dB_s)\over 1+r(x_s)},$$
 and let
\[\begin{array}{lll}
b_t&=&\half \sum_1^m\int_0^t \left( p(p-1){ [dr(X^i(x_s))]^2\over[1+ r(x_s)]^2}+ {p} {\nabla^2r(X^i(x_s), X^i(x_s))\over1+ r(x_s)}\right)ds\\
&+& p\int_0^t{dr(A^X(x_s))\over 1+r(x_s)}ds.
\end{array}\]

\noindent
We have:
\begin{eqnarray*}
[1+r(x_t)]^p=[1+r(x_0)]^p{\cal E}(M_t)\exp^{b_t}.\\
\end{eqnarray*}

\noindent
Here ${\cal E}(M_t)=\exp^{M_t-\half<M,M>_t}$. But $b_t$ is bounded from the
assumptions. So
$$E[1+r(x_t)]^p\le [1+ r(x_0)]^p\exp^{k_0[1+p^2]t}$$

\noindent
for some constant $k_0$. Thus
\begin{eqnarray*}
& &\sup_{x\in K} E\left( \exp^{6p^2 \int_0^t c\left[1+\ln(1+r(F_s(x))\right]ds}
\right)
\le {1\over t}\sup_{x\in K} \int_0^t E\left(\exp^{6p^2ct
\left[1+\ln(1+r(F_s(x))\right]}\right)ds \\
&\le& {1\over t}\exp^{6cp^2t}\sup_{x\in K} \int_0^t
 E\left([1+r(F_s(x))]^{6cp^2s}\right)ds<\infty
\end{eqnarray*}

So theorem ~\ref{th: existence}  applies to the function 
$f(x)=c[1+\ln(1+r(x))]$ to get the strong completeness.
\hfill\rule{3mm}{3mm}

\bigskip

\noindent{\bf Remarks:} 

\noindent
(i) From the above calculations we also get, for each $p>0$: 
$$P\{T_n<t\}\le {1\over n^p}[1+r(x_0)]^p \exp^{k_0 [1+p^2]t}.$$

\noindent
(ii) Note  that if the sectional curvatures are nonpositive, then
$\nabla^2r(x)\ge 0$ and so $\nabla^2 r(x)\le \triangle r(x)$. If the Ricci
 curvature has lower bound $-L^2(r(-))$, where $L$ is as before. Then 
\cite{Greene-Wu}
\begin{equation}
\triangle  r(-)\le (n-1)L(r(-))\hbox{coth}(rL(r(-)).
\label{inequality distance}
\end{equation}
In this case the theorem holds without further assumptions on the
 sectional curvatures.

\bigskip

In general, let $g\colon M\to R$ be a $C^2$ function, then
\[\begin{array}{lll}
\exp^{g(x_t)}&=&\exp^{g(x_0)} +\int_0^t \exp^gdg (X(x_s)dB_s)
+\half \int_0^t\exp^g \sum_1^m\left[dg(X^i(x_s))\right]^2 ds\\
&& +\int_0^t\exp^g\,\left( dg(A^X(x_s)) +
 \half \sum_1^m\nabla (dg)(X^i(x_s),X^i(x_s))\right)ds.\\
\end{array}\]
By Gronwall's inequality   $E\exp^{g(x_t)}<\exp^{g(x_0)}\exp^{kt}$ if
$dg(X^i)$ is bounded  for each $i$ and $\sum_1^m\nabla dg(X^i,X^i) + dg(A^X)$
 is bounded above. Using $g(x)=(1+r(x))^\epsilon$,  a similar proof  to that
 of Corollary 6.3 gives:

\begin{proposition}
Let $M$ be a complete Riemannian manifold with a pole. 
Assume its  sectional curvature is bounded from below by 
 $-L^2(r(-))$.  Here $L$ is a nondecreasing function bigger or equal to
 $1$. Then the SDE (1.1) is  complete  if for some $\epsilon>0$:
\begin{enumerate}
\item
 $|X(x)|^2\le {c[1+r(x)]^{2-\epsilon}\over  L(r(x))
\hbox{coth}\left(r(x)L(r(x))\right)}$;
\hskip 4pt
$|\nabla X(x)|^2\le c[1+(r(x))]^\epsilon$;
\item
$dr(A^X(x))\le c[1+r(x)]^{2-\epsilon}$;

\noindent
It is strongly complete, if also
\item 
 $H_p(x)(v,v)\le c[1+(r(x))]^\epsilon|v|^2$, for some $p>0$.
\end{enumerate}
\end{proposition}

Note that  this relaxes the conditions on the derivatives, compared to
theorem \ref{th: pole} but imposes more stringent bounds on the 
coefficients.

\section{ Strong completeness of nondegenerate equations}

 In this section we shall assume  that the SDS considered is a Brownian 
motion with drift $Z$, i.e.
$X^*X=Id$, and  $Z=:A^X=\half \sum_1^m\nabla X^i(X^i)+A$.

\noindent
Recall that $R$ is the curvature tensor and  Ric is the Ricci curvature. Then
$$\sum_1^m \left<R(X^i,v)(X^i),v \right>  =-{\rm Ric}(v,v),$$

\noindent        giving
\begin{equation}\begin{array}{ll} \label{eq: H for BM}
H_p(x)(v,v)=&2<\nabla Z(v),v>_x-Ric_x(v,v)+\sum_1^m |\nabla X^i(v)|_x^2 \\
	&+(p-2)\sum_1^m {1\over |v|^2} <\nabla X^i(v),v>_x^2.
\end{array}\end{equation}

\begin{theorem} \label{th: strong completeness Hessian}
Let $M$ be a complete Riemannian manifold. Assume  $|\nabla X|$  is bounded
 and $\half{\rm Ric}(v,v)- <\nabla Z(v),v>\ge -c|v|^2$ for some constant $c$. 
Then the Brownian motion  with drift $Z$ is strongly complete if complete.
\end{theorem}

\noindent
{\bf proof:} This follows from theorem  ~\ref{th: existence} by taking
 $f\equiv 1$.  \hfill \rule{3mm}{3mm}

In particular  suppose  the drift is $\nabla h$ for a smooth
 function $h$.  Then we have strong completeness  if $|\nabla X|$ is bounded
 and if  $\half$ Ric-Hess(h) is bounded from below, since a h-Brownian
 motion is complete if $\half$ Ric-Hess(h) is bounded from below.
  See   \cite{Bakry86}.

\bigskip

Let $p$ be a point in $M$. Let $r(x)$ denote the Riemannian distance 
 between $p$  and $x$.  The results in the last section hold for
 h-Brownian motions without the assumption  that there is  a pole
 for the manifold. Let $c$ be a constant.

\begin{theorem}\label{th: h-Brownian}
Let $M$ be a complete Riemannian manifold. Assume the Ricci curvature is
bounded from below by $-c(1+r^2(x))$.  Here $c$ is a constant. Suppose 
 $dr(Z(x)) \le c[1+r(x)]$ outside the cut locus cut(p) of $p$, 
 then the Brownian motion  with drift $Z$ is 
complete. Furthermore let $p>1$, then
 $E\left[r(x_t)\right]^p\le \left[1+r(x_0)\right]^p\exp^{k_0(1+p^2)t}$
for some constant $k_0$. It is strongly complete  and
$$\sup_{x\in K}E\left(\sup_{s\le t}|T_xF_s|^p\right)<\infty$$
for each $t>0$ and  compact set $K$, if the  following  also holds:

\noindent
1). $|\nabla X(x)|^2\le c[1+\ln(1+r(x))]$,

\noindent
2). Ric$_x(v,v)-2<\nabla Z(v),v>_x\, \ge -c[1+\ln(1+r(x))]|v|^2$.

\end{theorem}

\noindent {\bf Proof:}

\noindent
The proof of theorem 7.1 works here, noticing the following two points:

A. The  Ito formula for $[1+r(x_t)]^p$ (in the proof of theorem 
 \ref{th: pole})  holds with a correction term $L_t^p$:
 \begin{eqnarray*}
[1+r(x_t)]^p&=&[1+r(x_0)]^p+p\int_0^t [1+r(x_s)]^{p-1}dr(X(x_s)dB_s)\\
&+&{p(p-1)\over 2}\sum_1^m \int_0^t [1+r(x_s)]^{p-2} [dr(X^i(x_s)]^2ds \\
&+&{p\over 2} \sum_1^m \int_0^t [1+r(x_s)]^{p-1} \Delta r(x_s)ds\\
 &+&p\int_0^t [1+r(x_s)]^{p-1}dr(A^X(x_s))ds -L_t^p. \end{eqnarray*}

\noindent
where $L_t^p\ge 0$ and $\Delta r$ and $dr$ are defined to be zero on cut(p).
 See \cite{cr.ke.ma}.

B.  When  $x$ does not  belong to the  cut-locus $C(p)$ of $p$, there
is the following estimate from \cite{Greene-Wu}(p.26 and (2.27)  on p.35):
 $$|\triangle r(x)|\le (n-1)\sqrt{cL(r(x))}\hbox{coth}\left(r(x)
\sqrt{cL(r(x))}\right),$$
  Note also $\sum_1^m\nabla^2(X^i(x), X^i(x))=\triangle r(x)$. 
 However the cut-locus $C(p)$ has measure  zero, and so the  Brownian motion
 spends zero amount of time on the cut-locus by Fubini's theorem,  since it
 has a density with  respect to $dx$ for $dx$ the Riemannian volume  measure. 
 So  the proof of theorem \ref{th: pole} follows through.
\hfill\rule{3mm}{3mm}

Note that this method could also applied to the case of the  Ricci curvature
 is bounded below by 
$-L^2(r(-))$, where $L$ is a  nondecreasing  function greater
or equal to 1, just as in theorem \ref{th: pole}.

\bigskip

\noindent
{\bf Gradient Brownian systems}

Let $f: M \to R^m$ be an isometric embedding. Let $X(\cdot)(e)
=\nabla<f(\cdot),e>$. Such systems are called  gradient Brownian systems.
 Let $\nu_x$ be the space of normal vectors to $M$ at $x$. There is the second
 fundamental form:
$$\alpha_x: T_xM\times T_xM \to \nu_x$$

\noindent
and the shape operator: $A_x: T_xM\times \nu_x \to T_xM$
related by $\left<\alpha_x(v_1,v_2), w\right> =\left<A_x(v_1,w), v_2\right>$.
 Let  $\{e_i\}$ be an orthonormal basis for $R^m$. If $Y(x): R^m\to \nu_x$ is 
the orthogonal projection, then \cite{ELbook} \cite{ELflour}
$$\nabla X^i(v)=A_x \left(v,Y(x)e_i\right).$$
Let $f_1,\dots f_n$ be an o.n.b. for $T_xM$. Consider $\alpha_x(v,\cdot)$ as a linear map from $T_xM$ to $\nu_x$. Denote by $|\alpha_x(v,\cdot)|_{H,S}$ the corresponding Hilbert Schmidt norm, and $|\cdot|_{\nu_x}$ the norm of a vector in $\nu_x$. Accordingly we have:
\begin{eqnarray*}
&&\sum_1^m |\nabla X^i(v)|^2_x
=\sum_{i=1}^m\sum_{j=1}^n \left<A_x(v, Y(x)e_i), f_j\right>^2
=\sum_{i=1}^m\sum_{j=1}^n \left<\alpha_x(v,f_j), Y(x)e_i\right>^2\\
&&= \sum_{j=1}^m |\alpha_x(v,f_j)|_{\nu_x}^2
=|\alpha_x(v,\cdot)|_{H,S}^2,
\end{eqnarray*}
and
$$\sum_1^m \left <\nabla X^i(v), v\right>_x^2=|\alpha_x(v,v)|_{\nu_x}^2.$$

\noindent
This gives
\begin{equation}\begin{array}{cl}
H_p(v,v)=& -{\rm Ric}(v,v)+2<\nabla Z(v),v>+|\alpha_x(v,\cdot)|_{H,S}^2\\
 &+ {(p-2) \over |v|^2}|\alpha_x(v,v)|^2_{\nu_x}.
\end{array}
\end{equation}
Further,  Gauss's theorem:
${\rm Ric}(v,v)
=\left<\alpha(v,v), \hbox{trace } \alpha\right> -|\alpha(v,\cdot)|^2_{H,S}$
gives
\begin{equation}\begin{array}{cl}
H_p(v,v)=& -<\alpha(v,v), \hbox{trace } \alpha> 
	+ 2|\alpha_x(v,\cdot)|_{H,S}^2\\
&+{1\over |v|^2} (p-2)|\alpha_x(v,v)|^2_{\nu_x}
	+2<\nabla Z(v), v>_x.
\end{array}
\end{equation}

\bigskip

Thus the completeness and strongly completeness of a gradient Brownian motion
 rely only on  bounds on the second fundamental form and  on the 
drift:

\begin{corollary} \label{co: gradient}
Let $M$ be a closed immersed submanifold of $R^m$ with its second fundamental
 form  $\alpha$  bounded by  $c[1+\ln(1+r(x))]^{\half}$. Then a gradient 
 Brownian motion on  $M$   with drift $Z$  is strongly complete if
$$dr(Z)\le c[1+r(x)],$$
and
$$<\nabla Z(v),v>_x \, \le c[1+\ln(1+r(x))]|v|^2.$$
 It has a flow of diffeomorphisms if also 
$|Z(x)|\le c[1+r(x)]$, and
$|\nabla Z(x)|\le c[1+\ln(1+r(x))]$.
\end{corollary}

\noindent{\bf Proof:}
The strong completeness is clear from  theorem \ref{th: h-Brownian}.
 The  diffeomorphism property  comes from the fact that for gradient Brownian 
systems \cite{ELflow}, 
$$\sum_1^m\nabla X^i(X^i)=0.$$
So the 'adjoint' equation  (\ref{eq: adjoint}) to (\ref{eq: basic}) is also a 
 gradient Brownian system (with drift $-Z$).
\hfill\rule{3mm}{3mm}

Let $Z=0$, we get the following useful corollary:

\begin{corollary} Let $M$ be a complete Riemannian manifold isometrically
immeresed in $R^m$ with its
second fundamental form bounded by $c[1+\ln(1+r(x))]^{\half}$. Then 
the gradient Brownian motion on it has a flow  of diffeomorphisms.
\end{corollary}
See  also Baxendale \cite{Baxendale80} for a  discussion of flows
on manifolds with second fundamental form  bounded and globally Lipschitz.

\bigskip

According  to theorem \ref{pr: dim1}, a SDS is strongly 1-complete if
it is complete  and if $H_1(x)(v,v)\le c|v|^2$. But for gradient Brownian
 systems, we can do better. 
Let ${\cal E}(M_t)=\exp^{M_t^1-{<M^1,M^1>_t\over 2}}$, where $M_t^1$
is as defined before theorem \ref{th: existence} and let
$f(x)=\sup_{|v|= 1}H_1(x)(v,v)$.
Then we have:
\begin{proposition}
 Let $M$ be a closed immersed submanifold of $R^n$. Then a  stochastically
complete gradient Brownian system is strongly 1-complete if
$$\sup_{x\in K}E\left(\exp^{\half \int_0^Tf(F_s(x))ds}\right)<\infty$$
 for all compact set $K$ and  bounded stopping times $T$.
\label{th: Girsanov}
\end{proposition}
\noindent{\bf Proof:}
We shall use the notations of theorem \ref{pr: dim1}. Let 
$$\tilde {B_t}=B_t-\int_0^t Y(x_s)^*\left(\alpha_{x_s}
({v_s\over |v_s|},{v_s\over |v_s|})\right)ds$$
and let $\tilde {x_t}$ and $\tilde{v_t}$ be the solutions to the stochastic
differential equation
\begin{equation}\label{eq: Girsanoved}
dx_t=X(x_t)\circ d\tilde{B_t}+A(x_t)dt
\end{equation}
and the stochastic covariant equation:
$$dv_t=\nabla X(v_t)\circ d\tilde{B_t}+\nabla A(v_t)dt$$
respectively.  For $x\in M$, choose an o.n.b. $\{e_1, \dots, e_m\}$ for $R^m$, 
such that $\{X(x)(e_i)\}_1^n$ is an o.n.b. for $T_xM$ and $X(x)(e_j)=0$ for 
$j>n$. Then it is clear that  $X(Y^*(v))=0$ for $v\in \nu_x$. So equation
 (\ref{eq: Girsanoved}) is the same as our original stochastic differential 
equation (1), and thus $\tilde {x_t}$ has the same distribution as $x_t$ and
has no explosion.
On the other hand, by formula (\ref{eq: vt as exponential}):
\[\begin{array}{lll}
E|v_{S_j^K}|\chi_{S_j^K<t}&=&
|v_0|E\left({\cal E}(M_{t\wedge S_j^K})\exp^{a^1_{t\wedge S_j^K}}
\chi_{S_j^K<t}\right)\\
&=&|v_0|E\left({\cal E}(M_t)
\exp^{a^1_{t\wedge S_j^K}}\chi_{S_j^K<t}\right)\\
\end{array}\]

\noindent
by the optional stopping theorem. But by the Girsanov-Cameron-Martin formula
(\cite{ELbook}, \cite{RE-Yor91}),
\[\begin{array}{lll}
&&E\left({\cal E}(M_t)\exp^{a^1_{t\wedge S_j^K}}\chi_{S_j^K<t}\right)
=E\exp^{\half \int_0^{t\wedge S_j^K} H_1(\tilde{x_s})
({\tilde{v_s}\over |\tilde{v_s}|},{\tilde{v_s}\over |\tilde{v_s}|})ds}
\chi_{S_j^K<t}.\\
&\le& E\left(\exp^{\half \int_0^{t\wedge S_j^K} f(\tilde{x_s})ds }\chi_{S_j^K<t}\right)
=E\left(\exp^{\half \int_0^{S_j^K} f({x_s})ds }\chi_{S_j^K<t}\right)<\infty.
\end{array}\]

\noindent
Thus $\liminf_{j\to \infty}\sup_{x\in K}E|T_xF_{t\wedge S_j^K}|
\chi_{S_j^K<t}<\infty$, and the strong 1-completeness follows.
\hfill\rule{3mm}{3mm}

\section{Application to differentiation of semigroups}

 Assume the derivative of the  solution flow of equation $(\ref{eq: basic})$
 has  first moment:  $E|T_xF_s\chi_{s<\xi(x)}|<\infty$. We may define a
 semigroup (formally) of linear operators $\delta P_t$\index{\delta P_t$} on
 bounded measurable 1-forms as follows: for $v\in T_xM$ and $\phi$ a 
1-form
\begin{equation}
(\delta P_t)\phi(v) =E\phi\left(T_xF_t(v)\right)\chi_{t<\xi(x)}.
\label{eq: two1}
\end{equation}

It is in fact an $L^p$ semigroup under suitable conditions on the derivative 
flow $TF_t$. On the other hand, $\delta P_t(df)$ is clearly the formal
derivative of $P_tf$,  which can be checked to be true if the SDS
 concerned is strongly 1-complete and if $TF_t$ satisfy an integrability 
condition (see below). By virtue of 
the introduction of strong 1-completeness we can improve a theorem of Elworthy 
\cite{ELbook}. The assumption that the SDS is strongly 1-complete is, on the
other hand, a natural assumption: first  $dP_tf=(\delta P_t)(df)$ for $f\in 
BC^1$ implies completeness (take $f\equiv 1$), and in fact
 $dP_tf=(\delta P_t)(df)$ for
$f\in C_K^\infty$ and  $E|T_xF_t|\chi_{t<\xi(x)}<\infty$ implies completeness
\cite{application}. Here $BC^1$ is the space of bounded functions with
bounded continuous first derivatives.    And also strong 
1-completeness follows from  completeness if for a complete Riemannian metric 
$\sup_{x\in K}E\sup_{s\le t}|T_xF_s|<\infty$ for all compact sets $K$ (theorem
\ref{th: dim1}). For applications
of results in this section, see \cite{application}, and \cite{EL-LI}.

\begin{theorem}  $\label{th: differentiate semigroups}$
Assume strong 1-completeness. Suppose the map $r\to E|T_{\sigma(r)}F_t|$ 
is continuous for $r$ small, for all smooth curves
$\sigma\colon [0,\ell]\to M$.
If $f$ is $BC^1$, then $P_tf$ is $C^1$ and 
 $$d(P_tf)(x)=\delta P_t(df)(x).$$
\end{theorem}

\noindent {\bf Proof:}   
Let $x\in M$, $v\in T_xM$. Take a geodesic curve $\sigma\colon  [0,\ell]\to M$
 starting from $x$ with velocity $v$ such that the image set is contained in a
 compact neighbourhood $K$ of $x$.
 By the strong 1-completeness, $F_t(\sigma(s))$ is a.s. differentiable with respect to $s$. So for almostly all $\omega$: 
$${f\left(F_t\left(\sigma(s), \omega\right)\right)
-f\left(F_t(x, \omega)\right) \over s} 
={1\over s} \int_0^s df
\left(T_{\sigma(r)}F_t\left(\dot \sigma(r),\omega\right) \right)dr.$$

By the strong 1-completeness we know  $T_{\sigma(r)}F_t(\dot \sigma(r))$ is continuous in $r$ for almost all $\omega$. Thus:
\begin{eqnarray*}
&&E\lim_{s\to 0} {1\over s} \int_0^s df
\left(T_{\sigma(r)}F_s(\dot \sigma(r), \omega)\right)dr.
 = E\lim_{s\to 0} {1\over s} \int_0^s df(T_{\sigma(r)}
F_t(\dot \sigma(r),\omega))dr\\
&&= Edf(TF_t(v)).
\end{eqnarray*}
On the other hand, $ lim_{s\to 0}{1\over s}\int_0^s E|T_{\sigma(r)}F_t|dr
=E|T_xF_t|$  if the map
$r\to E|T_{\sigma(r)} F_t|$ is continuous. But $|df\left(T_{\sigma(r)}
F_t(\dot \sigma (r))\right)|\le |df|_\infty|T_{\sigma(r)}F_t|$, so
$\lim_{s\to 0}$$EI_s$  $=$$E\lim_{s\to 0}I_s$ giving
 $Edf(T_xF_t(v))=d(P_tf)(v)$.
 \hfill \rule{3mm}{3mm}

Let $\sigma(0)=x_0$, the required  continuity of the map $r\colon \to
E|T_{\sigma(r)}F_t|$ can be assured by one of the following conditions:
(1) There is a constant $\delta>0$ such that:
$$\sup_{x\in K}E |T_xF_t|^{1+\delta} <\infty,$$ 
for a compact neighbourhood  $K$ of $x_0$. 
(2) $E\sup_{x\in K}|T_xF_t|<\infty$ for a compact set $K$ containing $x_0$.

\bigskip

\begin{corollary}
Let $M$ be a complete Riemannian manifold. Suppose  SDS  (1) is complete 
 and satisfies:
$$H_{1+\delta}(v,v)\le k|v|^2.$$
Then $d P_tf=\delta P_t(df)$ if both $f$ and $df$ are bounded. Here $H$ is
as defined in section 5.
\end{corollary}

\noindent{\bf Proof:} 
First the system is strongly 1-complete by the boundedness of $H_1$. 
 On the other hand,  formula (\ref{eq: vt as exponential})
in section 5  gives:
$E|T_xF_t|^{1+\delta} \le \exp^{{c(1+\delta)\over 2}t}$.

\end{document}